\documentclass[aap]{amsart}

\input{mlstylemodif.sty}
\mathtoolsset{showonlyrefs}

\newtheorem{theorem}{Theorem}[section]
\newtheorem{remark}[theorem]{Remark}

\newtheorem{definition}[theorem]{Definition}

\hypersetup{colorlinks,
citecolor=blue,
citebordercolor=blue,
linkcolor=blue,
urlcolor=blue
}

\title[MFG Beyond Homogeneity and Anonymity]{An Overview of Some Extensions of Mean Field Games beyond Perfect Homogeneity and Anonymity}

\author{Mathieu Lauri\`ere}\thanks{NYU-ECNU Institute of Mathematical Sciences, NYU Shanghai, Shanghai, 200126, People’s Republic of China, email: zc3151@nyu.edu, mathieu.lauriere@nyu.edu}

\date{}

\usepackage[margin=1.35in]{geometry}

\begin{document}

\maketitle

\begin{abstract}
The mean field games (MFG) paradigm was introduced to provide tractable approximations of games involving very large populations. The theory typically rests on two key assumptions: \emph{homogeneity}, meaning that all players share the same dynamics and cost functions, and \emph{anonymity}, meaning that each player interacts with others only through their empirical distribution. While these assumptions simplify the analysis, they can be restrictive for many applications. Fortunately, several extensions of the standard MFG framework that relax these assumptions have been developed in the literature.
The purpose of these notes is to offer a pedagogical introduction to such models. In particular, we discuss multi-population MFGs, graphon MFGs, major-minor MFGs, and Stackelberg MFGs, as well as variants involving cooperative players.
\end{abstract}


\section{Introduction}

\subsection{Background}

The term \emph{mean field games} (MFGs for short) was coined by Lasry and Lions in~\cite{lasry2007mean} to describe a paradigm that provides tractable approximations of finite-player games. The mean-field approximation can be rigorously justified when the number of players is sufficiently large and the players are homogeneous and anonymous. Although assuming perfect homogeneity and anonymity is a significant restriction for many applications, it can be relaxed while preserving many of the advantages of MFG theory. In particular, several extensions of the ``standard'' setting have been introduced in the literature. 

The purpose of these notes is to offer a pedagogical introduction to several of these extensions. The aim is to convey the main modeling ideas in an accessible way rather than delve into mathematical details, which can be found in the references provided at the end of each section. It is our hope that these notes help stimulate interest in these extensions among a broad community.

The rest of this paper is organized as follows. Section~\ref{sec:notations} introduces useful notations. Section~\ref{sec:mfg} reviews the standard framework of MFGs. Sections~\ref{sec:mpmfg} and~\ref{sec:gmfg} are devoted to multi-population MFGs and graphon MFGs, respectively. Section~\ref{sec:big-player} discusses models with a more influential player, namely major-minor MFGs and Stackelberg MFGs. Each of these sections begins with a formulation of the finite-player game, then presents the asymptotic model, continues with a linear--quadratic example, and concludes with bibliographic notes. While the previous sections focus on non-cooperative players, Section~\ref{sec:coop} reviews several models with cooperative players, including mean field control, mean field-type games, mean field control games, and models combining cooperative and non-cooperative behaviors.

\subsection{Notations} 
\label{sec:notations}

For any integer $n$, we use the notation $[n] = \{1,\dots,n\}$. 
We will denote by $d$, $\ell$ and $p$ the dimension of the state, the action and the noise, respectively. 
Let $\actions$ denote the action set, which is assumed to be a subset of $\RR^\ell$. 
We denote by $\cP_2(\RR^d)$ the set of probability measures on $\RR$ having a finite second moment. 
Let $T \in [0,+\infty)$ be a finite time horizon. 
For any $x \in \RR^d$, $\delta_x$ denotes the Dirac mass at $x$. 
For a random variable $X$, we denote by $\cL(X)$ the law of $X$.

Throughout these notes, we prioritize the conceptual framework rather than the generality. For this reason, we focus on a relatively simple model. It should be mentioned that several variants have been considered in the literature, while \emph{still fitting in the perfectly homogeneous and anonymous framework}. They allow dealing with more complex and more realistic models of dynamics and cost functions, but are not essentially different as far as the main assumptions of homogeneity and anonymity are concerned. See the bibliographic notes in Section~\ref{sec:mfg-bibnotes} for more details.

\section{Mean field games} 
\label{sec:mfg}

We start our presentation with a relatively ``standard'' setup of MFGs. To motivate the problem, we first present the finite player game. %

\subsection{Finite population game}
\label{sec:standard-finite}
We first review the standard setting. The finite-player game is characterized by a tuple:
\begin{equation}
\label{eq:mfg-tuple-N}
    \left(\mu_0, b, \sigma, f, g, N\right),
\end{equation}
where:
\begin{itemize}
    \item $\mu_0 \in \cP_2(\RR^d)$ is the \textit{initial distribution},
    \item $b: \RR^d \times \actions \times \cP_2(\RR^d) \to \RR^d$ is the \textit{drift function},
    \item $\sigma \in \RR^{d \times p}$ is the \textit{diffusion coefficient},
    \item $f: \RR^d \times \actions \times \cP_2(\RR^d) \to \RR$ is the \textit{running cost function},
    \item $g: \RR^d \times \cP_2(\RR^d) \to \RR$ is the \textit{terminal cost function},
    \item $N \in \NN$ is the \textit{number of players}. 
\end{itemize} 

The state of player $i \in [N]$ at time $t \in [0,T]$ is denoted by $X^i_t \in \RR^d$. The empirical distribution at time $t \in [0,T]$ is:
\[
	\mu^{N}_t = \frac{1}{N} \sum_{i=1}^{N} \delta_{X^{i}_t}.
\]
The state of player $i \in [N]$ follows the dynamics:
\[
	dX^{i}_t = b(X^{i}_t, \alpha^{i}_t, \mu^{N}_t) dt + \sigma dW^{i}_t, \quad X^{i}_0 \sim \mu_0,
\]
where the initial states $X^i_0$ are independent, and $W^{i}, i=1,\dots,N$ are independent $p$-dimensional Brownian motions which are also independent of the initial states. Notice that the dynamics are coupled through the dependence of $b$ on the empirical distribution.

Given the control profile $\alpha^{-i} = (\alpha^1,\dots,\alpha^{i-1},\alpha^{i+1}, \dots, \alpha^N)$ for the other players, player $i \in [N]$ chooses $\alpha^i$ in a set of admissible controls to minimize:
\begin{equation}
\label{eq:def-JN}
	J_{N}(\alpha^{i}, \alpha^{-i}) = \mathbb{E} \left[ \int_0^T f(X^{i}_t, \alpha^{i}_t, \mu^{N}_t) dt + g(X^{i}_T, \mu^{N}_T) \right].
\end{equation}
The notion of solution is that of Nash equilibrium. Intuitively, it is a configuration in which no player has any incentive to change her control unilaterally: by doing so, she cannot reduce her cost. Formally, the definition is as follows.
\begin{definition}[Nash equilibrium]
\label{def:standard-finite-Nash}
A Nash equilibrium is a control profile $\underline{\hat{\alpha}} = (\hat{\alpha}^{i})_{i=1,\dots,N}$ such that: 
    \[
    	J_{N}(\hat{\alpha}^{i}, \hat{\alpha}^{-i}) 
	\leq J_{N}(\alpha^{i}, \hat{\alpha}^{-i}), \qquad \forall \alpha^i, \forall i \in [N].
    \]
\end{definition}

For most games, explicit analytical solutions are not available, and in general, computing a Nash equilibrium is a highly challenging problem, as demonstrated by results in complexity theory~\cite{daskalakis2009complexity,ummels2011complexity,daskalakis2023complexity}. This difficulty has motivated researchers to develop approximate models that are more tractable and yield controls corresponding to approximate Nash equilibria, wherein each player can gain only marginally by deviating unilaterally. More precisely, we introduce:
\begin{definition}[Approximate Nash equilibrium]
\label{def:approx-finite-Nash}
Let $\epsilon \ge 0$. An $\epsilon-$Nash equilibrium is a control profile $\underline{\hat{\alpha}} = (\hat{\alpha}^{i})_{i=1,\dots,N}$ such that: 
    \[
    	J_{N}(\hat{\alpha}^{i}, \hat{\alpha}^{-i}) 
	\leq J_{N}(\alpha^{i}, \hat{\alpha}^{-i}) + \epsilon, \qquad \forall \alpha^i, \forall i \in [N].
    \]
\end{definition}
Notice that an $\epsilon$-Nash equilibrium with $\epsilon=0$ is a Nash equilibrium. In the sequel, we focus on mean field games, which provide approximate Nash equilibria for finite-player games while being more tractable than Nash equilibria.

\subsection{Key assumptions}

The above model is formulated so as to be able to easily pass to the limit when the number of players, $N$, goes to infinity. This approach is in the spirit of propagation of chaos in statistical physics~\cite{sznitman1991topics}, which explains the terminology of \emph{mean field} games.

To be able to pass to the limit, we started our presentation with a class model that satisfies two important structural assumptions:
\begin{itemize}
	\item \textbf{Homogeneity:} all the players have the same form of dynamics and cost: the functions $b, f, g$ and the constant $\sigma$ (which could also be a function in general) do not depend on the player's index (but they depend on the payer's state). 
	\item \textbf{Anonymity:} player $i$ does not ``know'' the index of other players, meaning that she interacts with other players only through the empirical distribution $\mu^{N}_t$. 
\end{itemize}
Interactions through the empirical distribution are sometimes referred to as \emph{weak interactions} and can be viewed as a form of \emph{symmetry} with respect to the other players' states: instead of viewing $b$ as a function of $(X^{i}_t, \alpha^{i}_t, \mu^{N}_t)$, we can view it as a function of $(X^{i}_t, \alpha^{i}_t, X^{1}_t, \dots, X^{N}_t)$ which is symmetric with respect of the last $N$ inputs: for every permutation $\tau$ over $[N]$, for every $x \in \RR^d$, $a \in \cA$, $(x^1,\dots,x^N) \in \RR^d$, $b(x, a, x^1,\dots,x^N) = b(x, a, x^{\tau(1)},\dots,x^{\tau(N)})$. The same applies to all the functions involved in the model (namely, $f$ and $g$). As a consequence, we can expect that there exists a Nash equilibrium in which all the players use the same control. This simplifies the search for an equilibrium control.

These assumptions are instrumental to prove that the asymptotic game with infinitely many players (presented in the next subsection and called the MFG), provides a good approximation of the finite-player game (presented above). 
See the end for the next subsection for more details.

\subsection{Asymptotic game}
\label{sec:mfg-model}
When the number of players grows to infinity, we expect, at least informally, a form of propagation of chaos to hold. As a consequence, each player's state is less and less affected by the state of a specific other player. Instead, each player interacts only with the population distribution, also called mean field. In the limit $N \to +\infty$, we can reasonably assume the players' dynamics become decoupled, which should simplify the analysis. This intuition is formalized by the notion of mean field game (MFG).

The MFG is defined as follows. Given a mean field $\mu = (\mu_t)_{t \in [0,T]}$, $\mu_t \in \cP_2(\RR^d)$, and a control $\alpha$, the dynamics of a representative player's state when using control $\alpha$ are:
\begin{equation}
\label{eq:SDE-MFG}
    dX_t = b(X_t, \alpha_t, \mu_t) dt + \sigma dW_t, \quad X_0 \sim \mu_0.
\end{equation}
To stress the dependence on the control and the mean field, we will use the notation $X^{\alpha,\mu}$ for the solution of this SDE. 

The objective of such a representative player is to choose $\alpha$ to minimize the cost:
\begin{equation}
\label{eq:Jalphamu-meanfield}
    J(\alpha, \mu) = \mathbb{E} \left[ \int_0^T f(X^\alpha_t, \alpha_t, \mu_t) dt + g(X^\alpha_T, \mu_T) \right].
\end{equation}
Note that $J^N$ defined in~\eqref{eq:def-JN} is a function of a control profile, while $J$ is a function of the individual control and the mean field.

Given a control $\alpha$ and a mean field $\mu$, we will denote by $\mu^{\alpha,\mu}$ the mean field generated by $\alpha$ with mean field $\mu$, i.e., $\mu^{\alpha,\mu}_t = \cL(X^{\alpha,\mu}_t)$ for all $t \in [0,T]$.

\begin{definition}[MFG equilibrium]
\label{def:MFG-eq}
An MFG equilibrium (or simply mean field equilibrium) is a pair $(\hat{\alpha}, \hat{\mu})$ such that:
\begin{enumerate}
    \item Optimality: $\hat{\alpha}$ is a best response against $\hat\mu$, i.e.,
    \begin{equation}
        J(\hat{\alpha}, \hat{\mu}) \leq J(\alpha, \hat{\mu}), \qquad \forall \alpha.
    \end{equation}
    
    \item Consistency: The mean field $\hat\mu$ is the one generated by $\hat\alpha$, i.e.,
    \begin{equation}
        \hat{\mu}_t = \cL(X^{\hat{\alpha}, \hat\mu}_t), \quad \forall t \in [0,T].
    \end{equation}
\end{enumerate}
\end{definition}
This definition means that an MFG equilibrium is a fixed point. Alternatively, we can phrase the problem as follows.
\begin{definition}[MFG equilibrium -- equivalent definition]
\label{def:MFG-eq-ctrl}
$\hat\alpha$ is an MFG equilibrium control if and only if: $\hat{\alpha}$ is a best response against $\mu^{\hat\alpha}$, i.e.,
    \begin{equation}
        J(\hat{\alpha}, \mu^{\hat\alpha}) \leq J(\alpha, \mu^{\hat\alpha}), \qquad \forall \alpha.
    \end{equation}
\end{definition}

Here, every player is interested purely in her own cost function, taking the behavior of the rest of the population as fixed. If, on the contrary, all the players cooperate to minimize a common cost function, the solution concept becomes that of a social optimum, which has been studied under the terminology of mean field (type) control; see Section~\ref{sec:mfc}.

A natural question is: \textit{How is this MFG model related to the finite-player game presented above?} Under suitable conditions, two types of results can be obtained:
\begin{itemize}
	\item the MFG equilibrium control $\hat\alpha$ provides an $\epsilon$-Nash equilibrium control in the finite-player game in the sense of Definition~\ref{def:approx-finite-Nash} (i.e., if all of the $N$ players use $\hat\alpha$, then each player can be at most $\epsilon$-better off by choosing another control), and $\epsilon$ goes to $0$ as $N \to +\infty$; see e.g. the monographs~\cite[Chapter 5]{MR3134900} and~ \cite[Chapter 6]{carmona2018probabilistic2} (Theorems~6.7 and~6.13 without and with common noise, respectively);
	\item the $N$-player game equilibrium converges to the MFG equilibrium when $N \to +\infty$; see e.g. the monographs~\cite{MR3343705} and~\cite[Chapter 6]{carmona2018probabilistic2} (Theorems~6.18 and 6.28 for open-loop and closed-loop equilibria, respectively).
\end{itemize}
Generally, the second type of result is considered more difficult to obtain (and it requires at least existence of a Nash equilibrium in the finite-player game, which is not the case for the first type of result). 
See also the references provided in Section~\ref{sec:mfg-bibnotes} for detailed analyses. 

For completeness, we briefly discuss numerical methods. Definition~\ref{def:MFG-eq} can be interpreted as a fixed-point problem. This perspective suggests a direct numerical approach to computing the MFG solution by alternating between two steps: updating the control (by computing the best response to the current mean field) and updating the mean field (by determining the population distribution flow induced by the most recent best response). Such fixed-point iterations typically converge when the underlying mapping is a strict contraction. However, ensuring contractivity can be challenging in many models.
Alternative strategies include averaging over iterations (for instance, through fictitious play or online mirror descent) or regularizing the best response using an entropic penalty in the cost function. These two steps can be formulated either in terms of partial differential equations (namely, a Hamilton–Jacobi–Bellman equation for the value function, from which the best response is derived, and a Kolmogorov–Fokker–Planck equation for the evolution of the mean field) or in terms of stochastic differential equations, involving a backward stochastic differential equation for the value function or its derivative; see, for example,~\cite[Section 3.3]{carmona2018probabilistic} and~\cite{cardaliaguet2020introduction}, respectively. Alternatively, one can attempt to directly solve the forward–backward system, as in~\cite{achdou2010mean}.
The design of a numerical method typically involves three key components: selecting a model to approximate the relevant functions (for example, using a tabular representation or a class of function approximators such as neural networks), choosing a scheme to discretize the continuous equations (for example, a finite-difference method), and determining an approach to update the model parameters (for example, through linear algebra operations or stochastic gradient descent).
A comprehensive review of numerical methods for MFGs is beyond the scope of this survey. For more detailed discussions, we refer the reader to~\cite{achdou2010mean,achdoulauriere-2020-mfg-numerical,lauriere2021numerical,hu2024recent} and the references therein.

\subsection{Example}
\label{sec:mfg-ex}

As an example, let us take a linear-quadratic Gaussian (LQG) model, i.e., a model in which the drift is linear in the individual state, the individual action and the first moment of the mean field, the running cost is quadratic in these quantities, the terminal cost is quadratic in the state and the mean, and the state distribution is Gaussian at all times. So the interactions are purely through the first moment of the distribution. The fact that the distribution is Gaussian is not crucial for the characterization of the optimal control through ordinary differential equations (ODEs) but it helps to compute the optimal value and the evolution of the distribution, which reduces to the evolution of its mean and its variance in this case.

We will use the notation $\bar\mu = \int x \mu(dx)$. We borrow the following example from~\cite[Chapter 6]{MR3134900} and for ease of presentation, we focus on the one-dimensional case, i.e., $d=k=p=1$. 
We take:
\begin{align*}
	f(x,\mu,\alpha)	&= \frac{1}{2} \left[ Q x^2 + \bar{Q} \left(x - S \bar\mu \right)^2 + C \alpha^2 \right]
	\\
	g(x,\mu) 	&= \frac{1}{2} \left[ Q_T x^2 + \bar{Q}_T \left(x - S_T \bar\mu \right)^2 \right]
	\\
	b(x,\mu,\alpha) 	&= Ax + \bar{A} \bar\mu + B\alpha \, ,
\end{align*}
where $Q, C, \bar{Q}, Q_T, \bar{Q}_T$ are non-negative constants, and $A, \bar{A}, S, S_T$ and $B$ are constants. Let $\nu = \tfrac{1}{2} \sigma^2$.  We consider that the initial distribution is the normal distribution $\mu_0 = \mathcal{N}(\bar x_0, \sigma_0^2)$ for some $\bar{x}_0 \in \RR$ and $\sigma_0>0$.

\begin{remark}
\label{rem:ex-interpretation}
The model defined above can be extended to include more terms, see e.g.~\cite{graber2016linear}, but this one can already capture several interesting features. For example:
\begin{itemize}
	\item if $A=\bar{A}=0$ and $B=1$, then the drift is exactly the control, meaning that the agents control in which direction they move, up to the noise;
	\item if $A=-1$ and $\bar{A}=-1$, then $Ax + \bar{A} \bar\mu = (\bar\mu - x)$, which implies that the drift has a mean-reverting component, reminiscent of the Ornstein-Uhlenbeck process;
	\item likewise, if $S=1$ (resp. $S_T=1$), then the running cost (resp. terminal cost) gives an incentive to each agent to move towards the mean of the population.
\end{itemize}  
\end{remark}

Under suitable conditions on these coefficients, the MFG for the above model has a unique solution $(\hat{\alpha}, \hat{\mu})$ which satisfies the following. The proof relies on dynamic programming and on a suitable ansatz for the value function of an infinitesimal player when the population is in the Nash equilibrium; see \textit{e.g.}~\cite[Chapter 6]{MR3134900} for more details. We look for this value function in the form: $U(t,x) = \frac{1}{2}p_t x^2 + r_t x + s_t$. The coefficients must satisfy the ODE system:
 \begin{subequations}
     \begin{empheq}[left=\empheqlbrace]{align*}
		\int_\RR \xi \hat{\mu}(t,\xi) d \xi &= z_t,
		\\
		\hat{\alpha}(t,x) &= -B(p_t x + r_t) / C,
		\\
		J(\hat{\alpha}, \hat{\mu}) &= \int_\RR U(0,\xi) \mu_0(\xi) d\xi = \frac{1}{2} p_0(\sigma_0^2 + \bar x_0^2) + r_0 \bar x_0 + s_0,
		\\
		U(t,x) &= \frac{1}{2}p_t x^2 + r_t x + s_t,
     \end{empheq}
   \end{subequations}
where $(z, p, r, s)$ solve the following system of ODEs:
    \begin{subequations}
     \begin{empheq}[left=\empheqlbrace\,\,]{align}
     \label{AMS-num-eq:LQ-ODE-z}
     \frac{d z}{dt} &= (A+\bar{A} - \tfrac{B^2 p_t}{C}) z_t - \tfrac{B^2 r_t}{C}, 			&& z_0 = \int_{\RR} \xi \mu_0(\xi) d\xi,
     \\
     \label{AMS-num-eq:LQ-ODE-P}
     -\frac{d p}{dt} &= 2Ap_t - \tfrac{|Bp_t|^2}{C} + Q + \bar{Q}, 				&& p_T = Q_T + \bar{Q}_T,
     \\
     \label{AMS-num-eq:LQ-ODE-r}
     - \frac{d r}{dt} &= (A - \tfrac{B^2p_t}{C} )r_t + (p_t \bar{A} - \bar{Q} S) z_t, 	&& r_T = - \bar{Q}_T S_T z_T,
     \\
     \label{AMS-num-eq:LQ-ODE-s}
     - \frac{d s}{dt} &= \nu p_t - \frac{1}{2} \tfrac{|B r_t|^2}{C} + r_t \bar{A} z_t + \frac{1}{2}  \bar{Q} |Sz_t|^2, && s_T = \frac{1}{2} \bar{Q}_T |S_Tz_T|^2.
     \end{empheq}
   \end{subequations} 
   Here, $z$ is the mean of the population's distribution whereas $r$ (together with $p$) characterizes the best response. The last equation, for $s$, can be solved explicitly provided $z,p$ and $r$ are known.
   The second equation, for $p$, is a Riccati equation which does not involve the other variables. Under suitable assumptions on the coefficients of the problem, it admits a unique positive (symmetric if $d>1$) solution. 
   The most difficult part comes from the first and the third equations, namely for $z$ and $r$, which form a coupled forward-backward system. 
    This is a key difficulty which is at the core of MFGs. 
    As mentioned at the end of Section~\ref{sec:mfg-model}, one approach consists in solving the forward and the backward equations alternatively, possibly with averaging as in fictitious play. Another approach amounts to solve both equations simultaneously. We refer e.g. to~\cite[Section 2]{lauriere2021numerical} for more details and numerical experiments.

We refer e.g. to~\cite[Section 3]{graber2016linear} for more details in an even more general setting with interactions through the action distribution and common noise.

\subsection{Bibliographic notes}
\label{sec:mfg-bibnotes}

As already mentioned, the term mean field games was coined by Lasry and Lions~\cite{lasry2006jeux,lasry2006jeux2,lasry2007mean}. The idea of approximate Nash equilibria through an asymptotic model was studied around the same time in~\cite{huang2003individual,huang2005nash,huang2006large}. The theory has been developed by P.-L. Lions in his lectures at Coll\`ege de France. Since the standard MFG framework is not the main focus of these notes, we refer the interested reader e.g. to the monographs~\cite{MR3134900,gomes2016regularity,carmona2018probabilistic,carmona2018probabilistic2} and the surveys~\cite{gomes2014mean,cardaliaguet2020introduction,achdoulauriere-2020-mfg-numerical} and the lecture notes~\cite{cardaliaguet2018short,lacker2018mean}. See e.g.~\cite{dayanikli2024machinetutorial} for a recent tutorial with applications in operations research.

For the connection between finite-player games and MFGs, in addition to the monographs~\cite[Chapter 5]{MR3134900} and~\cite[Chapter 6]{carmona2018probabilistic2} mentioned above, see also~\cite{lacker2020convergence,lacker2023closed} for convergence results of closed-loop equilibria without and with common noise, respectively. Beyond asymptotic convergence, finer results have been established, including a central limit theorem~\cite{delarue2019mastercentral} and large deviation principles~\cite{delarue2020masterlarge}. Furthermore,~\cite{lauriere2022convergence,possamai2025non,jackson2025quantitative} proved the convergence of finite-player Nash equilibria to mean field equilibria for games with interactions through the distribution of actions.

For more details on the notions of anonymity, see e.g.~\cite[Chapter 2]{sandholm2010population} and \cite[Section 2.5.2]{nisan2007algorithmic}. For the notion of symmetric games, see e.g.~\cite[Section 2.1]{cardaliaguet2018short}.

Let us say a few words about variants of the above MFG that have been considered in the literature. These variants include more complex models, \emph{still satisfying homogeneity and anonymity}.   
We focus on the continuous time, finite horizon problem but other settings include the infinite horizon discounted setting and the ergodic setting, see e.g.~\cite[Sections 2.1--2.3]{lasry2007mean}, \cite[Section 7]{MR3134900} or \cite[Sections 2.4.2 and 2.4.3]{lauriere2022learning} in the discrete time setting; see~\cite{cardaliaguet2012long} for a rigorous justification of the connection between finite time and ergodic MFGs. Furthermore, we focus on interactions through the distribution of states only but the models could also include interactions through the distribution of actions or the joint distribution of states and actions; this is sometimes referred to as extended MFGs or MFG of controls, see e.g.~\cite{gomes2014existence,cardaliaguet2018mean,achdou2020mfgcontrols,kobeissi2022classical}. 
    We consider that the volatility is constant for simplicity; in general, it could depend on the state, the control and the mean field. Furthermore, the dynamics could include jumps, which is particularly relevant for models in which the state or some components of the state evolves in a discrete space. 
    Last, we restrict our attention to models in the spirit of classical stochastic optimal control problems, but other models include optimal stopping or impulse control.

\section{Multi-population MFGs} 
\label{sec:mpmfg}

We now turn to a first extension of the standard setup, in which the total population comprises several sub-populations which are homogeneous. To alleviate the presentation we will drop the ``sub-'' and simply call them populations. All the players are non-cooperative (whether it is within their population or with other populations). The standard setup is a special case, when there is just one population.

\subsection{Finite player model}

We consider $K$ populations of players. The finite-player multi-population game with $K$ populations is characterized by a tuple:
\begin{equation}
\label{eq:tuple-mpmfg-N}
    \left(\mu_{0,k}, b_k, \sigma_k, f_k, g_k, N_k\right)_{k=1,\dots,K},
\end{equation}
where:
\begin{itemize}
    \item $\mu_{0,k} \in \cP_2(\RR^d)$ is the \textit{initial distribution} of population $k$,
    \item $b_k: \RR^d \times \actions \times (\cP_2(\RR^d))^K \to \RR^d$ is the \textit{drift function} for population $k$,
    \item $\sigma_k > 0$ is the \textit{diffusion coefficient} for population $k$,
    \item $f_k: \RR^d \times \actions \times (\cP_2(\RR^d))^K \to \RR$ is the \textit{running cost function},
    \item $g_k: \RR^d \times (\cP_2(\RR^d))^K \to \RR$ is the \textit{terminal cost function},
    \item $N_k \in \NN$ is the \textit{number of players} in population $k$. 
\end{itemize} 

The empirical distribution of population $k \in [K]$ is defined as:
\[
	\mu^{N_k}_t = \frac{1}{N_k} \sum_{i=1}^{N_k} \delta_{X^{k,i}_t}.
\]
Player $i \in [N_k]$ in population $k \in [K]$ follows the dynamics:
\[
	dX^{k,i}_t = b_k(X^{k,i}_t, \alpha^{k,i}_t, \mu^{N_1}_t, \dots, \mu^{N_K}_t) dt + \sigma_k dW^{k,i}_t, \quad X^{k,i}_0 \sim \mu_{0,k}.
\]
We assume that the initial conditions $X^{k,i}_0, i\in[N_k], k \in[K]$ are sampled independently of each other, and that $W^{k,i}, i\in[N_k], k \in[K]$ are independent of each other and independent of the initial conditions. 

Given the control profile $\underline\alpha^{-(k,i)} = (\alpha^{k',i'})_{(k',i') \neq (k,i)}$ for other players, player $i$ in population $k$ wants to choose her control $\alpha^{k,i}$ to minimize the cost:
\[
	J^k_{N_1,\dots,N_K}(\alpha^{k,i}, \underline\alpha^{-(k,i)}) = \mathbb{E} \left[ \int_0^T f_k(X^{k,i}_t, \alpha^{k,i}_t, \mu^{N_1}_t, \dots, \mu^{N_K}_t) dt + g_k(X^{k,i}_T, \mu^{N_1}_T, \dots, \mu^{N_K}_T) \right].
\]

\begin{remark}
	In this model, the \emph{homogeneity} assumption is broken by the fact that $b,\sigma,f$ and $g$ can depend on the player's population. There is however homogeneity among the players of the same population. The \emph{anonymity} assumption is broken because $b,f$ and $g$ may depend differently on the empirical distributions of different populations. One could consider multi-population models breaking only one of the two assumptions. 
\end{remark}

The solution concept is that of a Nash equilibrium, in which every player cares only about her individual cost.
\begin{definition}[Nash equilibrium]
A Nash equilibrium consists of strategies $(\hat{\alpha}^{k,i})$ such that:
    \[
    	J^k_{N_1,\dots,N_K}(\hat{\alpha}^{k,i}, \hat{\alpha}^{-i}) 
    	\leq 
    	J^k_{N_1,\dots,N_K}(\alpha^{k,i}, \hat{\alpha}^{-i}), \quad \forall \alpha^{k,i}.
    \]
\end{definition}

\begin{remark}
\label{rem:MPMFG-to-MFG-finite}
	The (single-population) game presented in Section~\ref{sec:standard-finite} can be viewed as a special case of the multi-population game when $K=1$. Conversely, one can try to recast this multi-population game as a single population game by letting the population index $k$ be part of the state. In other words, one can consider $[K] \times \RR^d$ as a state space and write down a dynamics in which the first component does not evolve. One issue is that we want to have exactly $N_k$ players in population $k$ and this is not easy to capture using the single population setting. However, this issue will vanish in the mean field limit. See Remark~\ref{rem:MPMFG-to-MFG}. 
\end{remark}

\subsection{Asymptotic model}

In the asymptotic model, there are $K$ populations and population $k$ represents a proportion say $p_k \in [0,1]$ of the total population, with $\sum_{k=1}^K p_k = 1$. Intuitively, this corresponds to a situation where, in the finite player game presented above, $N_k/N \to p_k$ as $N_k$ and $N$ go to infinity. The limiting model is called a multi-population MFG (MPMFG for short).

Given the vector $(\underline\mu_t)_{t \in [0,T]} = (\mu^1_t, \dots, \mu^K_t)_{t \in [0,T]}$ of distribution flows for all the populations, the state of a representative player in population $k$ who uses control $\alpha^k$ has the dynamics:
\begin{equation}
    dX^k_t = b_k(X^k_t, \alpha_t^k, \underline\mu_t) dt + \sigma_k dW^k_t, \quad X^k_0 \sim \mu^k_0.
\end{equation}
We assume that the initial conditions $X^k_0, k \in[K]$ are sampled independently of each other, and that $W^k, k \in[K]$ are independent of each other and independent of the initial conditions. Since we focus on one representative player per population, we make use of $K$ different Brownian motions.  
We will use the notation $X^{k,\underline\mu,\alpha}$ if we want to stress the dependence on $\underline\mu$ and $\alpha$.

The objective of an player in population $k$ is to minimize the cost:
\begin{equation}
    J^k(\alpha, \underline\mu) = \mathbb{E} \left[ \int_0^T f_k(X^{k,\underline\mu,\alpha}_t, \alpha_t^k, \underline\mu_t) dt + g_k(X^{k,\underline\mu,\alpha}_T, \underline\mu_T) \right].
\end{equation}

\begin{remark}
\label{rem:MPMFG-to-MFG}
The issue mentioned in Remark~\ref{rem:MPMFG-to-MFG-finite} that one faces when trying to recast the multi-population game as a standard game disappears in the mean field limit. Indeed, in the mean field game, we can view the MPMFG as a single-population MFG, with state space $[K] \times \RR^d$ and initial distribution $\tilde\mu_0$ such that the first marginal of $\tilde\mu_0$ gives the proportion $p_k$ for each $k \in [K]$, and $\tilde\mu_0(\cdot|k) = \mu_{0}^k$. 
\end{remark}

\begin{definition}[Multi-Population MFG Equilibrium]
A multi-population MFG equilibrium is a collection $(\hat{\alpha}^k, \hat{\underline\mu}^k)_{k=1,\dots,K}$ such that:
\begin{enumerate}
    \item Optimality: Given the mean field $\hat{\underline\mu} = (\hat{\mu}^1, \dots, \hat{\mu}^K)$, for every $k$, the control $\hat{\alpha}^k$ is optimal for a representative player of type $k$, meaning:
    \begin{equation}
        J^k(\hat{\alpha}^k, \hat{\underline\mu}) \leq J^k(\alpha^k, \hat{\underline\mu}), \qquad \forall \alpha^k, \forall k \in [K].
    \end{equation}
    
    \item Consistency: For each population, the laws of optimally controlled processes matches the mean fields:
    \begin{equation}
        \hat{\mu}^k_t = \cL(X^{k,\hat{\alpha}}_t), \quad \forall t \in [0,T], \qquad \forall k \in [K].
    \end{equation}
\end{enumerate}
\end{definition}
Here, every agent of every population is interested only in her own cost function, taking the behavior of all the other agents as fixed. Alternatively, if the players of each population cooperate in order to minimize a common cost, the solution concept becomes that of a Nash equilibrium between cooperative groups of mean-field type, which has been studied under terminology of mean field type games; see Section~\ref{sec:mftg} for more details.

\subsection{Example}
\label{sec:mpmfg-ex}

We revisit the LQG example presented in Section~\ref{sec:mfg-ex} in the standard MFG setting, extending it with multi-population interactions. With the notation $\underline\mu = (\mu^1,\dots,\mu^K)$, we define: for $k \in [K]$,  
\begin{align*}
	f^k(x,\underline\mu,\alpha)	&= \frac{1}{2} \left[ Q^k x^2 + \bar{Q}^k \left(x - S^{k} \sum_{\ell=1}^K w_{k,\ell} \bar\mu^\ell \right)^2 + C^k \alpha^2 \right]
	\\
	g^k(x,\underline\mu) 	&= \frac{1}{2} \left[ Q_T^k x^2 + \bar{Q}_T^k \left(x -  S_T^{k} \sum_{\ell=1}^K w_{k,\ell} \bar\mu^\ell \right)^2 \right]
	\\
	b^k(x,\underline\mu,\alpha) 	&= A^kx + \bar{A}^{k} \sum_{\ell=1}^K w_{k,\ell} \bar\mu^\ell + B^k\alpha \, ,
\end{align*}
where the coefficients are generalizations of the ones in Section~\ref{sec:mfg-ex} since they can depend on the population index. Here, $w_{k,\ell} \in \RR_+$ captures the effect of population $\ell$ on a player in population $k$. We consider that the initial distribution of population $k$ is the normal distribution $\mu_0^k = \mathcal{N}(\bar x_0^k, (\sigma_0^k)^2)$ for some $\bar{x}_0^k \in \RR$ and $\sigma_0^k>0$.

In line with Remark~\ref{rem:ex-interpretation}, if for instance $A^k=-1$ and $\bar{A}^k=1$, then the representative player in population $k$ is attracted towards the weighted mean $\sum_{\ell=1}^K w_{k,\ell} \bar\mu^\ell$ of the population, whereas in a standard MFG, they would be attracted towards the (uniform) mean: $\tfrac{1}{K}\sum_{\ell=1}^K \bar\mu^\ell$.

Under suitable conditions on these coefficients, the MPMFG for the above model has a unique solution $(\hat{\underline\alpha}, \hat{\underline\mu}) = ((\hat{\alpha}^k)_{k \in [K]}, (\hat{\mu}^k)_{k \in [K]})$ which satisfies the following:
 \begin{subequations}
     \begin{empheq}[left=\empheqlbrace]{align*}
		\int_\RR \xi \hat{\mu}^k(t,\xi) d \xi &= z_t^k,
		\\
		\hat{\alpha}^k(t,x) &= -B^k(p_t^k x + r_t^k) / C^k,
		\\
		J^k(\hat{\alpha}^k, \hat{\underline\mu}) &= \int_\RR U^k(0,\xi) \mu_0^k(\xi) d\xi = \frac{1}{2} p_0^k((\sigma_0^k)^2 + (\bar x_0^k)^2) + r_0^k \bar x_0^k + s_0^k,
		\\
		U^k(t,x) &= \frac{1}{2}p_t^k x^2 + r_t^k x + s_t^k,
     \end{empheq}
   \end{subequations}
where $(z, p, r, s)$ solve the following system of ODEs:
    \begin{subequations}
     \begin{empheq}[left=\empheqlbrace\,\,]{align}
     \frac{d z^k}{dt} &= (A^k - \tfrac{|B^k|^2 p_t^k}{C^k}) z_t^k + \bar{A}^{k} \sum_{\ell=1}^K w_{k,\ell} z_t^\ell - \frac{|B^k|^2 r_t^k}{C^k}, 			&& z_0^k = \int_{\RR} \xi \mu_0^k(\xi) d\xi,
     \\
     -\frac{d p^k}{dt} &= 2A^kp_t^k - \tfrac{|B^k p_t^k|^2}{C^k} + Q^k + \bar{Q}^k, 				&& p_T^k = Q_T^k + \bar{Q}_T^k,
     \\
     - \frac{d r^k}{dt} &= (A^k - \tfrac{|B^k|^2 p_t^k}{C^k} )r_t^k + (p_t^k \bar{A}^{k} - \bar{Q}^k S^{k}) \sum_{\ell=1}^K w_{k,\ell} z_t^\ell, 	&& r_T^k = - \bar{Q}_T^k S_T^{k} \sum_{\ell=1}^K w_{k,\ell} z_T^\ell,
     \\
     - \frac{d s^k}{dt} &= \nu^k p_t^k - \frac{1}{2} \tfrac{|B^k r_t^k|^2}{C^k} + r_t^k \bar{A}^{k} \sum_{\ell=1}^K w_{k,\ell} z_t^{\ell} + \frac{1}{2} \bar{Q}^k \left|S^{k} \sum_{\ell=1}^K w_{k,\ell}  z_t^\ell\right|^2, && s_T^k = \frac{1}{2} \bar{Q}_T^k \left|S_T^{k} \sum_{\ell=1}^K w_{k,\ell} z_T^\ell\right|^2.
     \end{empheq}
   \end{subequations} 
We see that, here again, the equation for $p^k$ can be solved first, independently of the other variables. Then the equations for $(z^k)_{k \in [K]}$ and $(r^k)_{k \in [K]}$ should be solved together since they are coupled. Finally, the equation for $s^k$ can be solved last.

\subsection{Bibliographic notes}

We refer e.g. to~\cite{huang2006large,feleqi2013derivationmulti,cirant2015multi,MR3882529} for an analytical approach and to~\cite[Section 7.1.1]{carmona2018probabilistic} and~\cite{fujii2022probabilistic} for a probabilistic formulation. In the context of reinforcement learning, multi-population MFGs have been studied e.g. by~\cite{subramanian2020multi}. For applications, see e.g.~\cite{lachapelle2011mean} and~\cite[Section 6.1]{achdoulauriere-2020-mfg-numerical} for pedestrian crowds (see also~\cite[Section 5.2]{MR3882529} for a discussion on the non-cooperative and cooperative models), \cite[Sioux Falls example]{cabannes2022solving} for traffic routing, \cite{machtalay2025computational}  for traffic flow with different types of vehicles, \cite{kang2021task} for route planning with collision-avoidance, \cite{achdou2017mean,cardaliaguet2017segregation,barilla2021mean,camilli2024network} for urban planning, \cite{banez2020belief,ren2024hierarchical} for opinion dynamics, \cite{perolat2022scaling} for a predator-prey model, \cite[Section 6]{lachapelle2016efficiency,fujii2022mean} for price formation, \cite{casgrain2020mean} for algorithmic trading with differing beliefs, or \cite{vu2025heterogenous} for macroeconomic models. An application of MPMFG to clustering in machine learning has been studied in~\cite{aquilanti2021mean}. A reinforcement learning method for an LQ model has been studied in~\cite{uz2023reinforcement}. A more general model with heterogeneities has been considered in~\cite{cont2025homogenization}.

\section{Graphon mean field games}
\label{sec:gmfg}

Next, we consider a further generalization of multi-population MFGs in which each player can have a their own type and interact with each other players in a non-anonymous way. The interactions are encoded by a graph which, in the infinite-population limit, is replaced by a graphon. Such games have been referred to as \emph{graphon mean field games} or simply \emph{graphon games}.

\subsection{Finite player model}

The finite-player game is characterized by a tuple:
\begin{equation}
\label{eq:graph-game-tuple}
    \left(\mu_0, b, \sigma, f, g, N, w\right),
\end{equation}
where:
\begin{itemize}
    \item $\mu^i_0 \in \cP_2(\RR^d)$ is the \textit{initial distribution} for player $i \in [N]$,
    \item $b^i: \RR^d \times \actions \times \cP_2(\RR^d) \to \RR^d$ is the \textit{drift function} for player $i \in [N]$, 
    \item $\sigma^i \in \RR^{d \times p}$ is the \textit{diffusion coefficient} for player $i \in [N]$,
    \item $f^i: \RR^d \times \actions \times \cP_2(\RR^d) \to \RR$ is the \textit{running cost function} for player $i \in [N]$,
    \item $g^i: \RR^d \times \cP_2(\RR^d) \to \RR$ is the \textit{terminal cost function} for player $i \in [N]$,
    \item $N \in \NN$ is the \textit{number of players},
    \item $w \in \RR^{N \times N}$ is the \emph{weight matrix} representing the weights of a weighted graph; $w(i,j)$ represents the strength of interaction between player $i$ and player $j$.
\end{itemize} 

In line with most of the literature on graphon games, we will assume that $w$ is symmetric, i.e., $w(i,j) = w(j,i)$. However, this does not mean that the game is symmetric (in the sense introduced above): indeed, player $i$ may interact in a different way with different players, as we will see below.

Each player perceives an empirical distribution which depends on the graph of interaction. This distribution can be different for different players. The empirical neighborhood mean field of player $i \in [N]$ at time $t \in [0,T]$ is:
\[
	\mu^{i,N}_t = \frac{1}{N} \sum_{j=1}^{N} w_{i,j} \delta_{X^{j}_t}.
\]
Notice that this is not a probability distribution in general, without extra conditions on $w_{i,j}$ or a renormalization. 
The state $X^{i}$ of player $i \in [N]$ follows the dynamics:
\[
	dX^{i}_t = b^i(X^{i}_t, \alpha^{i}_t, \mu^{i,N}_t) dt + \sigma^i dW^{i}_t, \quad X^{i}_0 \sim \mu^i_0,
\]
where, as in the classical case, the initial states $X^i_0$ are independent, and the $W^{i}, i=1,\dots,N$ are independent $p$-dimensional Brownian motions which are also independent of the initial states.

Player $i$ chooses $\alpha^i$ in a set of admissible controls to minimize:
\begin{equation}
	J^w_{N}(\alpha^{i}, \alpha^{-i}) = \mathbb{E} \left[ \int_0^T f^i(X^{i}_t, \alpha^{i}_t, \mu^{i,N}_t) dt + g^i(X^{i}_T, \mu^{i,N}_T) \right].
\end{equation}
We use the superscript $w$ to stress the dependence on the weight matrix $w$, which occurs through the empirical neighborhood distribution $\mu^{i,N}$. 

\begin{remark}
	In this model, the \emph{homogeneity} assumption is broken by the fact that $b,\sigma,f$ and $g$ can depend on the player's index. The \emph{anonymity} assumption is broken because $b,f$ and $g$ may depend differently on different player, due to the graph of interactions. One could consider models breaking only one of the two assumptions. For instance, if $b$ is a function of the index but  $w$ is constant, then only the homogeneity assumption is broken; conversely, if $b,\sigma,f$ and $g$ are constant with respect to the index but $w$ is non-constant, then only the anonymity assumption is broken. 
\end{remark}

The notion of Nash equilibrium is defined in the same way as Definition~\ref{def:standard-finite-Nash}, with the new definition of $J_N$.

\subsection{Asymptotic model}

To formulate the asymptotic model when the number $N$ of players goes to infinity, we will use the concept of graphon mean field game (GMFG). Intuitively, this corresponds to situations where the graph is quite dense in the sense that each node has sufficiently many edges that the ``proportion'' of nodes that are neighbors does not vanish when the population size increases. 

In the asymptotic formulation, the parameters defining the game are the same as~\eqref{eq:graph-game-tuple} except that the players are indexed in $I=[0,1]$ and the weight matrix $w$ is replaced by a graphon, defined as follows.
\begin{definition}
	A graphon is a symmetric Borel-measurable function, $\WW : I\times I \rightarrow [0,1]$.
\end{definition}
This notion can be extended to other codomains. 
Formally, a graphon mean field game is defined by the tuple:
\begin{equation}
    \left(\mu_0, b, \sigma, f, g, \WW\right),
\end{equation}
where:
\begin{itemize}
    \item $\mu^u_0, b^u, \sigma^u, f^u, g^u$ are as before but for $u \in I$,
    \item $\WW: I\times I \rightarrow [0,1]$ is a graphon.
\end{itemize} 
Let $X^u_t$ denote the state of player $u \in I$ at time $t \in [0,T]$. Player $u$ is influenced by the aggregate:
\begin{equation}
    \label{eq:graphon-mu-u}
	\mu^{u,\WW}_t = \int_I \WW(u,v) \EE[\delta_{X^{v}_t}] dv.
\end{equation}
Note that this is not necessarily a probability measure, unless some extra conditions are imposed on $\WW$. 

The dynamics of player $u$ is given by: 
\[
	dX^{u}_t = b^u(X^{u}_t, \alpha^{u}_t, \mu^{u,\WW}_t) dt + \sigma^u dW^{u}_t, \quad X^{u}_0 \sim \mu^u_0,
\]
where, as in the classical case, the initial states $X^i_0$ are independent, and the $W^{u}, u \in I$ are $p$-dimensional Brownian motions such that $W^u$ is independent of $X^u_0$. To stress the dependence on $\alpha^u$, $\WW$ and $\underline\mu$, we will use the notation $X^{u,\alpha^u,\WW,\underline\mu}$.

\begin{remark}
    Intuitively, it might be more natural to replace the definition of $\mu^{u,\WW}$ in~\eqref{eq:graphon-mu-u} by: 
    \begin{equation}
    \label{eq:graphon-mu-u-exact}
        \mu^{u,\WW}_t = \int_I \WW(u,v) \delta_{X^{v}_t} dv,
    \end{equation}
    since player $u$ is expected to interact with player $v$ and not with the law of player $v$. However, this would require $X^v$ to be measurable with respect to $v$. 
    As explained in~\cite[Section 3.7]{carmona2018probabilistic}, it is not possible to have at the \emph{same time} that the Brownian motions $(W^u)_{u \in I}$ are independent and Lebesgue-measurable with respect to $u$. It is possible to have measurability and a weaker form of independence called essentially pairwise independence provided one is willing to work with a different probability space, using the framework of Fubini extensions~\cite{duffie2007existence}. Then, with a form of the law of large numbers called the exact law of large numbers~\cite{sun2006exact}, one can show that~\eqref{eq:graphon-mu-u-exact} is equivalent to~\eqref{eq:graphon-mu-u}. This approach has been used e.g. in~\cite{carmona2022stochastic,aurell2022stochastic,aurell2022finite}. However, although perhaps less intuitive, the definition~\eqref{eq:graphon-mu-u} of the aggregate $\mu^{u,\WW}$ also leads to an asymptotic game which can be showed to provide an approximate Nash equilibrium of the finite-player game with graph-based interactions. This can be viewed as a consequence of propagation of chaos-like results showing that, in the limit, individual noises average out. This viewpoint has been used e.g. in~\cite{caines2019graphon} for graphon games, \cite{bayraktar2023graphon,bayraktar2023propagation} for graphon dynamics, including forward-backward SDEs, which can be used to characterize graphon game equilibria.
\end{remark}

Given the distribution flows of all other players, represented by $\underline\mu = (\mu^v)_{v \in I}$, player $u \in I$ chooses $\alpha^u$ in a set of admissible controls to minimize:
\begin{equation}
\label{eq:JW-graphon-mf}
	J^{u,\WW}(\alpha^{u}, \underline\mu) = \mathbb{E} \left[ \int_0^T f^u(X^{u}_t, \alpha^{u}_t, \mu^{u,\WW}_t) dt + g^u(X^{u}_T, \mu^{u,\WW}_T) \right].
\end{equation}

\begin{definition}[Graphon mean field game equilibrium]
A graphon mean field game (GMFG) equilibrium is a pair $(\hat{\underline\alpha}, \hat{\underline\mu}) = ((\hat{\alpha}^u)_{u \in I}, (\hat{\mu}^u)_{u \in I})$ such that:
\begin{enumerate}
    \item Optimality: $\hat{\alpha}^u$ is a best response against $\hat{\underline\mu}$, i.e.,
    \begin{equation}
        J^{u,\WW}(\hat\alpha^{u}, \hat{\underline\mu}) \leq J^{u,\WW}(\alpha, \hat{\underline\mu}), \qquad \forall \alpha.
    \end{equation}
    
    \item Consistency: The mean field $\hat{\underline\mu}$ is the one generated by $\hat{\underline\alpha}$, i.e.,
    \begin{equation}
        \hat{\mu}^{u}_t = \cL(X^{u,\hat\alpha^u,\WW,\underline{\hat\mu}}_t), \quad \forall t \in [0,T], u \in I.
    \end{equation}
\end{enumerate}
\end{definition}

\begin{remark}
\label{rem:GMFG-to-MPMFG}
MPMFG can be viewed as a special case of GMFG.
Let us take $w$ as in the multi-population model of Section~\ref{sec:mpmfg-ex} and define $\WW(u,v) = K w_{k,\ell}$ for every $u,v \in I$ and $k,\ell \in [K]$ such that $(u,v) \in \left[\tfrac{k-1}{K}, \tfrac{k}{K}\right) \times \left[\tfrac{\ell-1}{K}, \tfrac{\ell}{K}\right)$. Then,  
\[
	\mu^{u,\WW}_t 
	= \int_I \WW(u,v) \EE[\delta_{X^{v}_t}] dv
	= \sum_{\ell \in [K]} \int_{\left[\tfrac{\ell-1}{K}, \tfrac{\ell}{K}\right)} \WW(u,v) \EE[\delta_{X^{v}_t}] dv
	= \sum_{\ell \in [K]} w_{k,\ell} K \int_{\left[\tfrac{\ell-1}{K}, \tfrac{\ell}{K}\right)} \EE[\delta_{X^{v}_t}] dv.
\]
So the interactions are only through the aggregates $\int_{\left[(\ell-1)/K, \ell/K\right)} \EE[\delta_{X^{v}_t}] dv$, $\ell \in [K]$. Hence we can expect that the equilibrium controls and mean fields are constant with respect to $u \in I$ on each sub-interval of the form $\left[(\ell-1)/K, \ell/K\right)$. We then recover the multi-population model of Section~\ref{sec:mpmfg}. 

\end{remark}

\begin{remark}
    In line with Remark~\ref{rem:MPMFG-to-MFG}, one can reformulate the GMFG as a standard MFG via a state-label formulation. The new state of a representative player at time $t$ comprises $X^u_t$ and the label $u$. See e.g.~\cite{cui2022learning,lacker2023label}. However, doing so might impose more restrictions on the regularity of the coefficients with respect to the label than what is necessary if the label is treated in an ad-hoc way, separately from $X^u_t$.
\end{remark}

\subsection{Example}
\label{sec:gmfg-ex}

We revisit the LQG example presented in Section~\ref{sec:mpmfg-ex} in the multi-population MFG setting, extending it with graphon interactions. 
With the notation $\underline\mu = (\mu^u)_{u \in I}$ and $\bar\mu^u = \int_\RR x \mu^u(dx)$, we define: for $u \in I$, 
\begin{align*}
	f^u(x,\underline\mu,\alpha)	&= \frac{1}{2} \left[ Q^u x^2 + \bar{Q}^u \left(x - \int_I S^u \WW(u,v) \bar\mu^v  dv \right)^2 + C^u \alpha^2 \right]
	\\
	g^u(x,\underline\mu) 	&= \frac{1}{2} \left[ Q_T^u x^2 + \bar{Q}_T^u \left(x - \int_I S_T^u \WW(u,v) \bar\mu^v dv \right)^2 \right]
	\\
	b^u(x,\underline\mu,\alpha) 	&= A^ux + \int_I \bar{A}^{u} \WW(u,v) \bar\mu^v dv + B^u\alpha \, ,
\end{align*}
where the coefficients are generalizations of the ones in Section~\ref{sec:mpmfg-ex} since they depend on an index $u \in I$ taking continuous values. We consider that the initial distribution of player $u$ is the normal distribution $\mu_0^u = \mathcal{N}(\bar x_0^u, (\sigma_0^u)^2)$ for some $\bar{x}_0^u \in \RR$ and $\sigma_0^u>0$.

\begin{remark}
    Here, we observe that the multi-population MFG example given in Section~\ref{sec:mpmfg-ex} is a special case of this one, obtained if $w_{k,\ell} = w_{\ell,k}$ and $\WW$ is piecewise constant, taking values $C_{k,\ell}$, $k,\ell \in [K]$, on sub-intervals of $I=[0,1]$. We come back to this point below, in the LQ example.
\end{remark}

Under suitable conditions on these coefficients, the GMFG for the above model has a unique solution $(\hat{\underline\alpha}, \hat{\underline\mu}) = ((\hat{\alpha}^u)_{u\in I}, (\hat{\mu}^u)_{u\in I})$ which satisfies the following:
 \begin{subequations}
     \begin{empheq}[left=\empheqlbrace]{align*}
		\int_\RR \xi \hat{\mu}^u(t,\xi) d \xi &= z_t^u,
		\\
		\hat{\alpha}^u(t,x) &= -B^u(p_t^u x + r_t^u) / C^u,
		\\
		J^u(\hat{\alpha}^u, \hat{\underline\mu}) &= \int_\RR U^u(0,\xi) \mu_0^u(\xi) d\xi = \frac{1}{2} p_0^u((\sigma_0^u)^2 + (\bar x_0^u)^2) + r_0^u \bar x_0^u + s_0^u,
		\\
		U^u(t,x) &= \frac{1}{2}p_t^u x^2 + r_t^u x + s_t^u,
     \end{empheq}
   \end{subequations}
where $(z, p, r, s)$ solve the following system of ODEs:
    \begin{subequations}
     \begin{empheq}[left=\empheqlbrace\,\,]{align}
     \frac{d z^u}{dt} &= (A^u - \tfrac{|B^u|^2 p_t^u}{C^u}) z_t^u + \bar{A}^{u} [\WW z_t]_{u} - \frac{|B^u|^2 r_t^u}{C^u}, 			&& z_0^u = \int_{\RR} \xi \mu_0^u(\xi) d\xi,
     \\
     -\frac{d p^u}{dt} &= 2A^up_t^u - \tfrac{|B^u p_t^u|^2}{C^u} + Q^u + \bar{Q}^u, 				&& p_T^u = Q_T^u + \bar{Q}_T^u,
     \\
     - \frac{d r^u}{dt} &= (A^u - \tfrac{|B^u|^2 p_t^u}{C^u} )r_t^u +  (p_t^u \bar{A}^{u}  - \bar{Q}^u S^{u} ) [\WW z_t]_{u}, 	&& r_T^u = - \bar{Q}_T^u S_T^{u} [\WW z_T]_{u},
     \\
     - \frac{d s^u}{dt} &= \nu^u p_t^u - \frac{1}{2} \tfrac{|B^u r_t^u|^2}{C^u} + r_t^u \bar{A}^u [\WW z_t]_{u} + \frac{1}{2} \bar{Q}^u    \left|S^u [\WW z_t]_{u} \right|^2, && s_T^u = \frac{1}{2} \bar{Q}_T^u \left|S_T^u [\WW z_T]_{u}\right|^2,
     \end{empheq}
   \end{subequations} 
   with $[\WW z_t]_{u} = \int_I \WW(u,v) z_t^v dv$ for brevity. 
We see that, here again, the equations for $p^u$ can be solved first, then the equations for $z^u$ and $r^u$ are coupled, and finally the equations for $s^u$ can be solved. 

Taking $\WW$ as in Remark~\ref{rem:GMFG-to-MPMFG}, one can check that the solutions to the above ODEs are constant (with respect to $u \in I$) on each sub-interval of the form $\left[(\ell-1)/K, \ell/K\right)$. We then recover the multi-population ODE system of Section~\ref{sec:mpmfg-ex}.

\subsection{Bibliographic notes}

Delarue studied in~\cite{delarue2017mean} a finite-player game with Erd{\"o}s-Renyi graph, which converges to a standard mean field game. The term graphon game was coined by Parise and Ozdaglar in~\cite{parise2019graphon,parise2023graphon}. The framework has attracted a lot of interest in the engineering community, see e.g.~\cite{caines2021graphon,gao2021linear,tchuendom2020master,foguen2022optimal,gao2023lqg,caines2023critical,foguen2024infinite}.  For the connection between finite-player game and graphon game using Fubini extensions, see e.g.~\cite{carmona2022stochastic,aurell2022stochastic}. For the convergence of finite-particle systems towards graphon systems, see e.g.~\cite{bayraktar2023graphon,bayraktar2023propagation,bayraktar2024concentration}.
Graphon dynamics and graphon games with jumps have been studied in~\cite{amini2025graphon,amini2023stochastic}. For applications, see e.g.~\cite{tangpi2024optimal} in finance, \cite{aurell2022finite} for epidemic management, or \cite{liu2025modeling} for rumor propagation. An LQ graphon game with common noise has been studied in~\cite{xu2024linear}.  On the numerical side, see~\cite{aurell2022finite,lauriere2025deep} for deep learning methods in which neural networks are used to learn functions of the index, without discretizing the interval $I=[0,1]$. Learning and reinforcement learning methods have been studied in~\cite{cui2022learning,zhang2023learning,zhou2024graphon,dong2025last}. As mentioned above, graphon games correspond to finite-player games with graph-based interactions when the graph is relatively dense. Extensions include hypergraphon games~\cite{cui2022hypergraphon}, directed graphon games~\cite{fabian2022meandigraphon}, or games with sparser interaction graphs, see e.g.~\cite{lacker2022case,fabian2023learning,crucianelli2024interacting,fabian2024learninggraphex}.

\section{More influential players: Major-Minor and Stackelberg MFG}
\label{sec:big-player}

Recall that in MFGs (as well as in MPMFGs and GMFGs), every player is infinitesimal and has no influence on the mean field: they perform their optimization taking the mean field as given. In this section, we present an extension of MFGs in which there is a player whose influence on the rest of the population is not negligible.

We will describe the general structure of the dynamics and the cost functions. Then we will distinguish between two notions of solutions: Nash equilibrium (aka major-minor MFG) and Stackelberg equilibrium (aka Stackelberg MFG).

\begin{remark}
    To alleviate the presentation, we will consider the population of infinitesimal players to be homogeneous and anonymous, as in a standard MFG of Section~\ref{sec:standard-finite}; but one could combine the presence of an influential player with e.g. GMFGs. 
    Also, for the sake of simplicity, our presentation will stick to the case of a single influential player, but the models could include several such influential players, with a suitable notion of solution (e.g., Nash equilibrium) between them.  
\end{remark}

\subsection{Finite player model}

We consider a system with $N$ identical players and one influential player. We will often use superscript $0$ to denote quantities related to this player. The model is characterized by a tuple:
\begin{equation}
\label{eq:tuple-mfg-influential-player}
    \left(\mu_0, b, \sigma, f, g, N, \mu_0^0, b^0, \sigma^0, f^0, g^0\right),
\end{equation}
where the first $6$ components are for the population and the last $5$ are for the influential player:
\begin{itemize}
    \item $\mu_0 \in \cP_2(\RR^d)$ is the \textit{initial distribution},
    \item $b: \RR^d \times \actions  \times \cP_2(\RR^d) \times \RR^{d^0} \times \actions^0 \to \RR^d$ is the \textit{drift function},
    \item $\sigma \in \RR^{d \times p}$ is the \textit{diffusion coefficient},
    \item $f: \RR^d \times \actions \times \cP_2(\RR^d) \times \RR^{d^0} \times \actions^0 \to \RR$ is the \textit{running cost function},
    \item $g: \RR^d \times \cP_2(\RR^d) \times \RR^{d^0} \times \actions^0 \to \RR$ is the \textit{terminal cost function},
    \item $N \in \NN$ is the \textit{number of players},
    \item $\mu_0^0 \in \cP_2(\RR^{d^0})$ is the \textit{influential player's initial distribution},
    \item $b^0: \RR^{d^0} \times \actions^0  \times \cP_2(\RR^d)  \to \RR^d$ is the \textit{influential player's drift function},
    \item $\sigma^0 \in \RR^{d^0 \times p^0}$ is the \textit{influential player's diffusion coefficient},
    \item $f^0: \RR^{d^0} \times \actions^0\times \cP_2(\RR^d)  \to \RR$ is the \textit{influential player's running cost function},
    \item $g^0: \RR^{d^0} \times \cP_2(\RR^d) \times \actions^0 \to \RR$ is the \textit{influential player's terminal cost function}.
\end{itemize} 

The state of player $i \in [N]$ at time $t \in [0,T]$ is denoted by $X^i_t \in \RR^d$, and the state of the influential player at time $t \in [0,T]$ is denoted by $X^0_t \in \RR^{d^0}$. The empirical distribution at time $t \in [0,T]$ is:
\[
	\mu^{N}_t = \frac{1}{N} \sum_{i=1}^{N} \delta_{X^{i}_t}.
\]
The state of player $i \in [N]$ follows the dynamics:
\[
	dX^{i}_t = b(X^{i}_t, \alpha^{i}_t, \mu^{N}_t, X^{0}_t, \alpha^{0}_t) dt + \sigma dW^{i}_t, \quad X^{i}_0 \sim \mu_0,
\]
and the state of the influential player follows the dynamics:
\[
	dX^{0}_t = b^0(X^{0}_t, \alpha^{0}_t, \mu^{N}_t) dt + \sigma^0 dW^{0}_t, \quad X^{0}_0 \sim \mu_0^0,
\]
where the initial states $(X^i_0)_{0,1,\dots,N}$ are independent, $W^{0}$ is a $p^0$-dimensional Brownian motion, and the $W^{i}, i=1,\dots,N$ are $p$-dimensional Brownian motions; $W^0,W^1,\dots,W^N$ are assumed to be independent and independent of the initial states. Notice that the influential player's dynamics depends on the $N$-player population only through its empirical distribution $\mu^N$.

Given the control profile $\alpha^{-i} = (\alpha^1,\dots,\alpha^{i-1},\alpha^{i+1}, \dots, \alpha^N)$ for the other players and the control $\alpha^0$ of the major player, player $i \in [N]$ chooses $\alpha^i$ in a set of admissible controls to minimize:
\begin{equation}
	J_{N}(\alpha^{i}, \alpha^{-i}, \alpha^0) = \mathbb{E} \left[ \int_0^T f(X^{i}_t, \alpha^{i}_t, \mu^{N}_t, X^{0}_t, \alpha^{0}_t) dt + g(X^{i}_T, \mu^{N}_T, X^0_T, \alpha^0_T) \right],
\end{equation}
and given the control profile $\underline\alpha = (\alpha^i)_{i\in[N]}$ for the players, the influential player chooses $\alpha^0$ to minimize:
\begin{equation}
	J_{N}^0(\alpha^{0}, \underline\alpha) = \mathbb{E} \left[ \int_0^T f^0(X^{0}_t, \alpha^{0}_t, \mu^{N}_t) dt + g^0(X^{0}_T, \mu^{N}_T, \alpha^0_T) \right].
\end{equation}

\begin{remark}
	In this model, the \emph{homogeneity} assumption is broken by the fact that $b,\sigma,f$ and $g$ can depend on the player's index because they are different for the influential player (with index $0$ in our notations). The \emph{anonymity} assumption is broken because $b,f$ and $g$ may depend differently on this influential player. However, all the players in the population are treated in the same way. If we remove the influential player, then we obtain the MFG of Section~\ref{sec:standard-finite}. Instead of a homogeneous and anonymous game, one could also consider multi-population or graphon games discussed above and add an influential player. 
\end{remark}

We will consider two types of solution concepts. Before that, let us introduce the corresponding mean field model.

\subsection{Asymptotic model}

The asymptotic model is defined by the same tuple~\eqref{eq:tuple-mfg-influential-player} except that $N$ is not needed anymore. In contrast with standard MFG, see Section~\ref{sec:mfg-model}, we cannot fix a mean field and define the influential player's cost, because this player will have an influence on the population: when the influential player optimizes her control, she wants to take into account the impact it will have on the population. As a consequence, we will define the cost functions as functions of controls. 

Consider a control $\alpha$ for a representative player of the population, a control $\alpha^{MF}$ used by the population, and a control $\alpha^0$ used by the influential player. Then, we consider the following dynamics, where $X, X^{MF}$ and $X^0$ represent the state of the representative player using control $\alpha$, the state of a representative player in the population using control $\alpha^{MF}$, and the state of the influential player:
\begin{equation}
    \begin{cases}
        dX_t = b(X_t, \alpha_t, \mu_t, X^{0}_t, \alpha^{0}_t) dt + \sigma dW_t, \quad X_0 \sim \mu_0,
        \\
        dX^{MF}_t = b(X^{MF}_t, \alpha^{MF}_t, \mu_t, X^{0}_t, \alpha^{0}_t) dt + \sigma dW_t, \quad X^{MF}_0 \sim \mu_0,
        \\
        dX^{0}_t = b^0(X^{0}_t, \alpha^{0}_t, \mu_t) dt + \sigma^0 dW^{0}_t, \quad X^{0}_0 \sim \mu_0^0,
    \end{cases}
\end{equation}
where $\mu_t = \cL(X^{MF}_t|W^0)$, which is the conditional law of the state $X^{MF}$ given the noise $W^0$ affecting the influential player. Notice that the three dynamics involve $\mu_t$. Moreover the dynamics of $X$ depends on $\mu_t$ and $X^{0}$, but the converse is not true: one can solve for $(X^{MF}, X^0)$ first and then find $X$. To stress the dependence on the controls used in the dynamics, we will use the notation: $X^{\alpha,\alpha^{MF}, \alpha^0}$, $X^{MF,\alpha^{MF},\alpha^0}$, and $X^{0,\alpha^{MF},\alpha^0}$. We will also denote $\mu_t^{\alpha^{MF},\alpha^0} = \cL(X_t^{MF,\alpha^{MF},\alpha^0}|W^0)$.

Given controls $\alpha^{MF}$ and $\alpha^0$ for the population and the influential player respectively, the objective of a representative player in the population is to choose $\alpha$ to minimize the cost:
\begin{align}
    J(\alpha, \alpha^{MF}, \alpha^0) 
    &= \mathbb{E} \Big[ \int_0^T f(X^{\alpha,\alpha^{MF}, \alpha^0}_t, \alpha_t, \mu^{\alpha^{MF}, \alpha^0}_t, X^{0,\alpha^{MF}, \alpha^0}_t, \alpha^{0}_t) dt 
    \\
    &\qquad\qquad + g(X^{\alpha,\alpha^{MF}, \alpha^0}_T, \mu^{\alpha^{MF}, \alpha^0}_T, X^{0,\alpha^{MF}, \alpha^0}_T, \alpha^0_T) \Big].
\end{align}
Given a control $\alpha^{MF}$ for the population, the objective of the influential player is to choose $\alpha^0$ to minimize the cost:
\begin{align}
    J^0(\alpha^0, \alpha^{MF}) 
    &= \mathbb{E} \left[ \int_0^T f^0(X^{0,\alpha^{MF}, \alpha^0}_t, \alpha^{0}_t, \mu^{\alpha^{MF}, \alpha^0}_t) dt + g^0(X^{0,\alpha^{MF}, \alpha^0}_T, \mu^{\alpha^{MF}, \alpha^0}_T, \alpha^0_T) \right].
\end{align}
As mentioned above, the fact that we define the influential player's cost as a function of $\alpha^{MF}$ and not as a function of the mean field directly is not insignificant: it implies that if $\alpha^0$ changes while $\alpha^{MF}$ is fixed, the mean field might still change. Notice that this would still be true if only $X^0$ and not $\alpha^0$ was appearing in the drift $b$ of the population.

\subsection{Solution concepts}

Now that we have defined the dynamics and the cost functions for the population of players and for the influential player, there are several notions of solutions that can be studied. Here we will focus on two notions. Each notion can be relevant in different applications. In the first two subsections, we proceed formally without discussing precisely the class of controls. Different types of information structures are discussed in Subsection~\ref{sec:info-structure}.

\subsubsection{Major-minor MFG}

One possibility is to consider that the influential player and all the players in the mean field population are playing a Nash equilibrium. In this case, the influential player is often referred to as a \emph{major player} and the other players are referred to as \emph{minor players}.

The solution concept is defined as follows in the finite-player model.
\begin{definition}
A Nash equilibrium for the finite-player major-minor game is a control profile $\underline{\hat{\alpha}} = (\hat{\alpha}^{i})_{i=1,\dots,N}$ for the minor players and a control $\hat\alpha^0$ for the major player such that: 
    \[
    	J_{N}(\hat{\alpha}^{i}, \hat{\alpha}^{-i}, \hat\alpha^0) 
	    \leq J_{N}(\alpha^{i}, \hat{\alpha}^{-i}, \hat\alpha^0), \qquad \forall \alpha^i, \forall i \in [N],
    \]
and
    \[
    	J_{N}^0(\hat{\alpha}^{0}, \hat{\underline\alpha}) 
	    \leq J_{N}^0(\alpha^{0}, \hat{\underline\alpha}), \qquad \forall \alpha^0.
    \]
\end{definition}

In the asymptotic model, the notion of Nash equilibrium is defined as follows, which is closer to Definition~\ref{def:MFG-eq-ctrl} than Definition~\ref{def:MFG-eq}, except for the fact that the costs are defined in terms of the population's control instead of the mean field, for reasons discussed above.
\begin{definition}[MMMFG Equilibrium]
A major-minor MFG (MMMFG or M3FG) equilibrium is a pair $(\hat{\alpha}, \hat{\alpha}^0)$ such that:
\begin{enumerate}
    \item Minor player optimality: $\hat{\alpha}$ is a best response for the representative minor player against $(\hat\alpha,\hat\alpha^0)$, i.e.,
    \begin{equation}
        J(\hat{\alpha}, \hat\alpha,\hat\alpha^0) \leq J(\alpha, \hat\alpha,\hat\alpha^0), \qquad \forall \alpha.
    \end{equation}
    
    \item Major player optimality: $\hat{\alpha}^0$ is a best response for the major player against $\hat\alpha$, i.e.,
    \begin{equation}
        J^0(\hat{\alpha}^0, \hat\alpha) \leq J^0(\alpha^0, \hat\alpha), \qquad \forall \alpha^0.
    \end{equation}
\end{enumerate}
\end{definition}

\subsubsection{Stackelberg MFG}
A different point of view consists in interpreting the influential player as a leader who can choose a control first and, conditioned on this, the population reacts and plays a Nash equilibrium. We will call the influential player the \emph{leader} and the players in the population the \emph{followers}.

In the finite-player game, the solution is defined as follows. It will be convenient to use the following notation: for a control $\alpha^0$ for the leader, $\mathcal{NE}_N(\alpha^0)$ denotes the set of $N$-player Nash equilibria for the followers, i.e., control profiles $\underline{\hat{\alpha}} = (\hat{\alpha}^{i})_{i=1,\dots,N}$ such that:
    \[
    	J_{N}(\hat{\alpha}^{i}, \hat{\alpha}^{-i}, \alpha^0) 
	    \leq J_{N}(\alpha^{i}, \hat{\alpha}^{-i}, \alpha^0), \qquad \forall \alpha^i, \forall i \in [N].
    \]
\begin{definition}
A Stackelberg equilibrium (or leader-follower equilibrium) for the finite-player game is a control $\alpha^{*,0}$ for the major player and a profile $\underline{\hat{\alpha}} = (\hat{\alpha}^{i})_{i=1,\dots,N}$ for the population such that:
\begin{itemize}
    \item Optimality for the leader: $\alpha^{*,0}$ is optimal for her when the population reacts by playing a Nash equilibrium:
    \[
    	\min_{\hat{\underline\alpha} \in \mathcal{NE}_N(\alpha^{*,0})} J_{N}^0(\alpha^{*,0}, \hat{\underline\alpha}) 
	    \leq \min_{\hat{\underline\alpha} \in \mathcal{NE}_N(\alpha^0)} J_{N}^0(\alpha^{0}, \hat{\underline\alpha}), \qquad \forall \alpha^0.
    \]
    \item Nash equilibrium for the followers: $\underline{\hat{\alpha}} \in \argmin_{\hat{\underline\alpha} \in \mathcal{NE}_N(\alpha^{*,0})} J_{N}^0(\alpha^{*,0}, \hat{\underline\alpha})$.
\end{itemize}
\end{definition}
Note that, above, the leader's goal is optimistic in the sense that she tries to optimize her control assuming that the population will play the best possible Nash equilibrium from her point of view (i.e., the one with the minimal cost). It would be possible to consider other definitions, such as a pessimistic viewpoint, in which the leader tries to optimize her cost with respect to the worst possible Nash equilibrium (this amounts to replace the $\min$ by a $\max$ in the definition).

The asymptotic game is sometimes referred to as Stackelberg MFG (SMFG) or leader-follower MFG. To define it, we first introduce the notation: for a control $\alpha^0$ for the leader, $\mathcal{NE}(\alpha^0)$ denotes the set of MFG equilibria for the followers, i.e., controls $\hat{\alpha}$ such that:
\[
    J(\hat{\alpha}, \hat{\alpha}, \hat\alpha^0) 
    \leq J(\alpha, \hat{\alpha}, \hat\alpha^0), \qquad \forall \alpha.
\]
Now, we introduce the SMFG solution concept.
\begin{definition}
A Stackelberg MFG equilibrium (or leader-follower MFG equilibrium) for the asymptotic game is a control $\alpha^{*,0}$ for the major player and a control $\hat{\alpha}$ for the followers:
\begin{itemize}
    \item Optimality for the leader: $\alpha^{*,0}$ is optimal for her when the population reacts by playing a Nash equilibrium:
    \[
    	\min_{\hat{\alpha} \in \mathcal{NE}(\alpha^{*,0})} J^0(\alpha^{*,0}, \hat{\alpha}) 
	    \leq \min_{\hat{\alpha} \in \mathcal{NE}(\alpha^0)} J^0(\alpha^{0}, \hat{\alpha}), \qquad \forall \alpha^0.
    \]
    \item Nash equilibrium for the followers: $\hat{\alpha} \in \argmin_{\hat{\alpha} \in \mathcal{NE}(\alpha^{*,0})} J^0(\alpha^{*,0}, \hat{\alpha})$.
\end{itemize}
\end{definition}

\subsection{Information structure}
\label{sec:info-structure}

In contrast with standard MFGs, where open-loop and closed-loop (Markovian) controls are usually equivalent in the sense that they lead to the same equilibria, the presence of a more influential player raises questions about the information available to each player.

Let us start with major–minor MFGs. One may consider, for instance, an open-loop control setting where the representative minor player’s control (denoted by $\alpha$) depends on the filtration generated by her own noise process ($W$) as well as the major player’s noise process ($W^0$), while the major player’s control ($\alpha^0$) depends only on the filtration generated by her own noise process (and not on the minor player’s noise, as it is negligible in the mean-field limit). This is, for example, the approach followed in~\cite{carmona2016probabilistic} and~\cite[Section 7.1]{carmona2018probabilistic2}. \cite{carmona2017alternative} defined major–minor equilibria as a fixed-point problem in the space of controls and, by searching for controls of the form $\alpha_t = \alpha(t, W_{[0,T]}, W^0_{[0,T]})$ and $\alpha^0_t = \alpha^0(t, W^0_{[0,T]})$, recovered, in an LQ model, the Nash equilibrium found in~\cite{carmona2016probabilistic}. 
Alternatively, one could consider “Markovian” closed-loop controls where, for instance, the representative minor player’s control is a function of her individual state, the major player’s state, and the (stochastic) mean field, while the major player’s control is a function of her individual state and the mean field, i.e., $\alpha_t = \alpha(t, X_t, X^0_t, \mu_t)$ and $\alpha^0_t = \alpha^0(t, X^0_t, \mu_t)$. This approach is adopted, for example, in~\cite[Section 2.3]{carmona2017alternative} and~\cite{cardaliaguet2020remarks}. More generally, the controls could depend on the full state trajectories, i.e., $\alpha_t = \alpha(t, X_{[0,t]}, X^0_{[0,t]}, \mu_t)$ and $\alpha^0_t = \alpha^0(t, X^0_{[0,t]}, \mu_t)$, as in the general formulation proposed in~\cite[Section 2.2]{carmona2017alternative}. It should be noted that solving major–minor MFG models is much more involved than solving standard MFGs, even in LQ settings; see, for example,~\cite[Section 7.1.6]{carmona2018probabilistic2} for a specific case and its resolution via a forward–backward system of SDEs and ODEs. Regarding numerical methods, \cite{carmona2017alternative} compute the solution to an LQ model by solving the associated system of ODEs. In discrete-time, finite-state settings, \cite{cui2024learning} propose an algorithm based on an adaptation of fictitious play. To the best of our knowledge, there is still no numerical method available for general continuous-time, continuous-space major–minor MFGs.

For Stackelberg MFGs, following~\cite{elie2019tale}, which draws inspiration from the literature on contract theory, part of the literature has focused on situations where the leader has no state and her control appears only in the terminal cost of the followers, that is, only $\alpha^0_T$ plays a role (in other words, the drift and running cost of the followers are constant with respect to $\alpha^0_t$). In such cases, the terminal action is interpreted as a terminal payment given to the agent based on her trajectory, i.e., $\alpha^0_T = \alpha^0(X_{[0,T]})$. A generalization allows the leader to influence the followers at any time, meaning that $\alpha^0_t$ affects the followers’ drift and running cost; see, for example,~\cite{dayanikli2024machine}, which also proposes a deep learning method for such problems. In a specific class of models where the followers’ terminal cost coincides with the leader’s terminal payment, \cite{elie2019tale} show that the Stackelberg equilibrium reduces to a mean field control problem (see Section~\ref{sec:mfc} below for more background on this class of problems). The authors also provide several examples with explicit or semi-explicit solutions. However, in general, formulating optimality conditions for Stackelberg MFGs is much more challenging than for standard MFGs, and further research is needed to develop numerical methods and establish their convergence properties.

\subsection{Bibliographic notes}

For more background on MFGs with major and minor players, see for instance~\cite[Section 7.1]{carmona2018probabilistic}. Such models have been considered, e.g., by \cite{huang2010large, nguyen2012linear, wang2012distributed20122093, bensoussan2016mean,firoozi2020convex,firoozi2020epsilon} for LQG systems, by \cite{nourian2013mm, sen2016mean} in non-linear and partially observed settings, and more recently by \cite{carmona2016probabilistic, carmona2017alternative, lasry2018mean, cardaliaguet2020remarks}. The continuous-time, finite-state setting has been studied in~\cite{carmona2016finite}. Extensions include risk-sensitive models~\cite{chen2023risk,liu2023lqg}, unobserved latent processes~\cite{firoozi2018mean}, recursive utilities~\cite{huang2024mean} and impulse control~\cite{chen2024deciding}. \cite{cui2024learning} has proposed a learning algorithm based on fictitious play. For applications, see e.g.~\cite{firoozi2017execution} in finance, \cite{feron2020price} in electricity markets, \cite{dayanikli2023multi} in electricity production, \cite{lin2022optimal} in electric vehicles charging. Recently, \cite{delarue2025major} has shown that adding a form of common noise helps to ensure existence, uniqueness and stability of Nash equilibria in major-minor MFGs. A model combining leader-follower and major-minor structures has been proposed in~\cite{sanjari2025incentive}.

In the spirit of contract theory, MFG models involving a leader and a large population of players have been studied in~\cite{elie2019tale} and~\cite{carmona2021finite}, respectively in continuous state space and finite state space. In the language of contract theory, the leader is called principal and the players are called agents. In the continuous space setting, several works have focused on LQ models, see e.g.~\cite{nourian2012mean,MR3376121,bensoussan2017linearstackelberg,moon2018linearstackelberg,lin2018openstackelberg,oner2018mean,yang2021linear,wang2024leader,chen2025sampled,wang2025linear}. Several extensions of the framework discussed in these notes have also been considered, such as partial observability~\cite{si2025linear,si2024linearcommon}, leader with backward dynamics~\cite{si2021backward,cong2024direct}, delays~\cite{MR3376121}, terminal state constraint~\cite{fu2018meanfield}, and informational uncertainties~\cite{xiang2024stochastic}. Different game structures have also been considered. For instance \cite{hu2022principal} considered a model with a single agent and a mean field of principals, \cite{oner2018mean} studied several layers of leader-follower structures, and \cite{moller2016constrained} considered a Stackelberg game between two infinite populations of non-cooperative players.
Various applications of Stackelberg MFG have been considered in the literature. In particular, let us mention advertisement \cite{salhab2016,salhab2018dynamic,carmona2021mean,salhab2022dynamic}, %
epidemic control \cite{aurell2022optimal}, energy demand response~\cite{elie2021mean,mastrolia2025agency}, renewable energy certificate markets~\cite{firoozi2021principal}, carbon taxes~\cite{carmona2022meancarbon,bichuch2024stackelberg},  heating, ventilation, and air conditioning~\cite{zhang2025stackelberg}, and macro-economic models~\cite{nuno2017bank,nuno2017optimal}. 
As in classical MFG, one can expect that the mean field model is close to the corresponding finite-player model (with one leader and a finite but large population of players); this connection has been studied rigorously e.g. in~\cite{wang2014hierarchical6716957}, \cite{bergault2023mean} and~\cite{djete2023stackelberg} using respective analytical and probabilistic techniques. 
\cite{guo2022optimizationstackelberg} proposed an optimization viewpoint for Stackelberg MFG in finite spaces. 
To the best of our knowledge, there are very few works on numerical methods for Stackelberg MFGs; we refer to~\cite{campbell2021deep} for a machine learning method in a specific setting, and to~\cite{dayanikli2024machine} for a deep learning approach tackling more general models.

\section{Cooperative models}
\label{sec:coop}

In this section, we review how the aforementioned models can be adapted to include cooperative interactions. We follow the structure adopted in the previous sections and we briefly mention each model without going into much details.

\subsection{From MFG to mean field control}
\label{sec:mfc}

The standard MFG model described in Section~\ref{sec:mfg}. The finite-player model is still described by the same tuple~\eqref{eq:mfg-tuple-N}, whose components have the same interpretation. What changes is the notion of solution: Instead of a Nash equilibrium (Definition~\ref{def:standard-finite-Nash}), we consider a \emph{social optimum}. To this end, we introduce the average cost over the population:
\[
    J_{N}^{SO}(\underline\alpha) = \frac{1}{N} \sum_{i=1}^N J_{N}(\alpha^{i}, \alpha^{-i}).
\]

\begin{definition}[Social optimum]
\label{def:standard-finite-Social}
A social optimum is a control profile $\underline{\alpha}^* = (\alpha^{*i})_{i=1,\dots,N}$ such that: 
    \[
    	J_{N}^{SO}(\underline\alpha^{*}) 
	\leq J_{N}^{SO}(\underline\alpha), \qquad \forall \underline\alpha = (\alpha^{i})_{i=1,\dots,N}.
    \]
\end{definition}
In other words, a control profile is a social optimum if it is optimal for the whole population, compared with any other control profile (not just the ones obtained with unilateral deviations).

Letting the number of players go to infinity, we obtain a mean field problem which has been called mean field type control, mean field control (MFC), or control of McKean-Vlasov dynamics. We recall that $J$ is defined in~\eqref{eq:Jalphamu-meanfield} as a function of one player's control and the mean field. In MFC, all the players are assumed to use the same control. So we now introduce the following cost, which is a function of the control only:
\begin{equation}
    J^{SO}(\alpha) = J(\alpha,\mu^\alpha) 
    = \mathbb{E} \left[ \int_0^T f(X^\alpha_t, \alpha_t, \mu_t^\alpha) dt + g(X^\alpha_T, \mu_T^\alpha) \right],
\end{equation}
where $\mu^\alpha$ is the mean field obtained when all the players use the control $\alpha$, i.e., $\mu^\alpha_t = \cL(X^{\alpha}_t)$, where $X^\alpha$ solves:
\begin{equation}
    dX_t = b(X_t, \alpha_t, \mu^\alpha_t) dt + \sigma dW_t, \quad X_0 \sim \mu_0.
\end{equation}
Note that, differently from the SDE used in MFG, namely, \eqref{eq:SDE-MFG}, here the mean field is not fixed: it is determined endogeneously since it depends on the law of the SDE solution itself. This is hence a McKean-Vlasov (MKV) equation. 

The MFC problem is then defined as follows. 
\begin{definition}[MFC optimum]
\label{def:MFC-opt}
An MFC optimum is a control $\alpha^*$ which is optimal for $J^{SO}$. 
\end{definition}

Under suitable conditions, it can be shown that the MFC optimal control provides an approximately optimal control for the finite-player cooperative problem.

The ratio of the cost obtained by an average player in MFG versus MFC is called the price of anarchy, by analogy with the literature on classical games~\cite{roughgarden2005selfish}; see~\cite{carmona2019price} in the context of MFGs. The comparison between MFG and MFC has been highlighted in several other works, including~\cite{carmona2013control}, \cite{cardaliaguet2019efficiency} using the notion of inefficiency, and \cite{wang2019mean} for a model of production output control with sticky prices. \cite{carmona2023nash} studied how selfish players can be incentivized to behave as in a mean-field social optimum, and how, in the absence of such incentives, they deviate from a social optimum to adopt a Nash equilibrium.

For more background on MFC, see e.g. the monographs~\cite[Chapter 4]{MR3134900}, \cite[Chapter 6]{carmona2018probabilistic} For numerical aspects, see e.g.~\cite{achdoulauriere-2020-mfg-numerical,lauriere2021numerical}. In particular, for numerical illustrations of the price of anarchy, see~\cite[Sections 2.6 and 4.4]{lauriere2021numerical}, in an LQ model and a model of crowd motion.

\subsection{From multi-population MFG to mean field type games and mixed population models}
\label{sec:mftg}

Next, we turn to models with several sub-groups. We present three different solution concepts. 
\begin{enumerate}
    \item Just as we did for multi-population MFG in Section~\ref{sec:mpmfg}, it is of course possible to consider an \emph{MFC problem with several populations}. In this case, there are several groups defined by the same tuple as~\eqref{eq:tuple-mpmfg-N} but the solution concept is that of a global social optimum, in which all the players of all the populations cooperate to minimize a common average cost. See e.g. \cite{xie2022social,feng2023unified} in LQ models, \cite[Section 2]{MR3882529} for a PDE approach and~\cite{fujii2022probabilistic} for a probabilistic approach. 
    \item Alternatively, we can consider that the the players inside each group are cooperative with each other and try to minimize the average cost of the group, but the groups do not cooperate with each other. This leads to the notion of \emph{mean field type games} (MFTGs) studied in~\cite{tembine2017mean}. In the asymptotic model where each group has infinitely many players (but there is a finite of group), it means that each group solves an MFC problem but between the groups one looks for a Nash equilibrium. Applications include blockchain token economics~\cite{barreiro2019blockchain}, risk-sensitive control~\cite{tembine2015risk} or, more generally, engineering~\cite{djehiche2017mean,barreiro2021meanbook}. Similar problems have also been called mean field games among teams~\cite{subramanian2023mean}, team-against-team mean field problems~\cite{sanjari2023nash,yuksel2024information}, and teamwise mean field competitions~\cite{yu2021teamwise}. The case of zero-sum MFTG has received special interest, see~\cite{MR3882529,basar2021zero,cosso2019zero,guan2024zero} for the theoretical foundations and~\cite{carmona2020policy,carmona2021linear,zaman2024independent} for numerical aspects using policy gradient. A connection between robust MFC and zero-sum MFTGs has been presented in~\cite{zaman2024robust}; see also~\cite{xu2025robust} for robust MFC. But the framework of MFTGs also covers general sum games with more than two (mean-field) coalitions; see e.g.~\cite{shao2024reinforcement} for finite-state models and reinforcement learning methods. The above models are for a finite number of groups. In the limit where the number of groups goes to infinity, MFTGs lead to an MFG model with infinitely many players where each player solves an MFC problem. Such models have been studied under the terminology of (mixed) \emph{mean field control games} in~\cite{angiuli2022reinforcementmfcgicaif-no-award,angiuli2023reinforcement,angiuli2024analysis}, motivated by multi-scale reinforcement learning algorithms. Similar games have also been studied in~\cite{li2022dynamic}, in~\cite[Section 3.2.2]{carmona2023nash} in a special case corresponding to an interpolation between MFG and MFC, and in~\cite{dayanikli2024how,dayanikli2025cooperation} under the terminology of \emph{mixed individual mean field model}.
    \item Last, we can also imagine a mixed model in which some groups are composed of cooperative players and some groups are non-cooperative. This has been studied under terminology of \emph{mixed-population mean field model} in~\cite{dayanikli2024how,dayanikli2025cooperation}, with applications to the tragedy of the commons.
\end{enumerate}

\subsection{From graphon MFG to graphon mean field control}

Next, the graphon-based model discussed in Section~\ref{sec:gmfg} can also be studied from the point of view of a cooperative solution concept, i.e., a social optimum. In this case, the goal is to find a control profile $\underline\alpha^* = (\alpha^{*u})_{u \in I}$ which minimizes the following social cost, integrated over the population:
\begin{equation}
	J^{SO,\WW}(\underline\alpha) = \int_I J^{u,\WW}(\alpha^u, \underline\mu^{\underline\alpha}) du,
\end{equation}
where $\underline\mu^{\underline\alpha}$ is the mean field generated by the control profile $\underline\alpha$, and  we recall that $J^{u,\WW}$ is defined in~\eqref{eq:JW-graphon-mf}, as a function of one player's control and the mean field.

To the best of our knowledge, social optimum problems in graphon-type interactions have received much less attention than Nash equilibria in graphon games. Such problems have been studied under the terminology of graphon mean field control (GMFC), or optimal control for non-exchangeable systems.

In the LQ setting, see~\cite{gao2018graphon,gao2019optimal,liang2021finite,liang2024asymptotically,xu2024social,pham2025lqgmfc}. \cite{dunyak2024quadratic} considered a different model than the one described in these notes, with a Q-noise. GMFC problems were studied in \cite{hu2023graphon} in discrete time and finite state spaces. \cite{pham2024gmfc} studied GMFC problems via PDE methods, analyzing the associated Bellman dynamic programming equations, while \cite{cao2025probabilistic} developed the probabilistic analysis of GMFC and a Pontryagin maximum principle, as well as the a propagation-of-chaos type result.

\subsection{Influential player with MFC population} 

Last, the models of Section~\ref{sec:big-player} with more influential players can also be adapted to the case with a cooperative mean-field population. 

Major-minor MFC models have been considered e.g. in~\cite{cui2024major}. Stackelberg equilibria with an MFC population have been studied for instance in~\cite{hubert2022incentives} for epidemic management, in~\cite{huang2021social,feng2024stackelberg,cong2024linear} and~\cite[Section 5.1]{jiang2025brief} for LQ models. Several extensions have been considered, see for instance ~\cite{lv2023linear} with regime switching, \cite{zhang2024stackelberg} in the context of $H_\infty$ control, and \cite{li2024linear} with partial observation.

\section{Conclusion}

In these notes, we have presented an overview of several extensions of the standard MFG and MFC frameworks that move beyond the assumptions of perfect homogeneity and anonymity. The references listed at the end of each section, while not exhaustive, illustrate the significant progress achieved in each of these directions.

Several avenues remain open for future research. First, the theoretical foundations are still incomplete for certain extensions, such as major–minor MFGs and Stackelberg MFGs in general settings. Second, numerical methods for these models have received less attention compared to those for standard MFGs and MFCs; major–minor and Stackelberg formulations, in particular, present substantial computational challenges. Third, these extensions are primarily motivated by practical applications: relaxing the assumptions of the classical mean-field framework offers the potential to construct more realistic models. Although some applications have already appeared in the literature, developing more sophisticated models and applying them to real-world problems remains an exciting direction for future work.

\section*{Acknowledgements}
This work was sparked by questions from the audience at a workshop organized by Tielong Shen at Dalian University of Technology. The author is grateful for the valuable questions and comments. The work has also greatly benefitted from the comments of the anonymous reviewers.

\bibliographystyle{abbrv}
\bibliography{referencesall}

\end{document}